\documentclass{article}
\usepackage{latexsym}
\usepackage{graphicx}

\begin{document}
\begin{center}
\textbf{A NEW FRACTIONAL DERIVATIVE WITHOUT SINGULAR KERNEL: APPLICATION TO
THE MODELLING OF THE STEADY HEAT FLOW}
\end{center}

\begin{center}
by
\end{center}

\begin{center}
\textit{Xiao-Jun YANG}$^{a}$, \textit{H. M. SRIVASTAVA
}$^{b,c}$, \textit{J. A. Tenreiro MACHADO}$^{d}$
\end{center}

\begin{center}
a Department of Mathematics and Mechanics, China University of Mining and
\end{center}

\begin{center}
Technology, Xuzhou 221008, People's Republic of China
\end{center}

\begin{center}
b Department of Mathematics and Statistics, University of Victoria,
Victoria, British Columbia V8W 3R4, Canada
\end{center}

\begin{center}
c China Medical University, Taichung 40402, Taiwan, ROC
\end{center}

\begin{center}
d Department of Electrical Engineering, Institute of Engineering,
Polytechnic of Porto, Rua Dr. Antonio Bernardino de Almeida, Porto 4249-015,
Portugal
\end{center}

\textit{In this article we propose a new fractional derivative without singular kernel. We consider the potential application for modeling the steady heat-conduction problem. The analytical solution of the fractional-order heat flow is also obtained by means of the Laplace transform. }


\textbf{Introduction}

Fractional derivatives with singular kernel [1], namely, the
Riemann-Liouville [2-3], Caputo [4-5] and other derivatives (see [6-8] and
the references therein), have nowadays a wide application in the field of
heat-transfer engineering.

More recently, the fractional Caputo-Fabrizio derivative operator without
singular kernel was given as (see[1,9-12]):
\begin{equation}
\label{eq1}
{ }^{CF}D_x^{\left( \nu \right)} T\left( x \right)=\frac{\left( {2-\nu }
\right)\Im \left( \nu \right)}{2\left( {1-\nu } \right)}\int\limits_0^x
{\exp \left( {-\frac{\nu }{1-\nu }\left( {x-\lambda } \right)}
\right)T^{\left( 1 \right)}\left( \lambda \right)d\lambda } ,
\end{equation}
where $\Im \left( \nu \right)$ is a normalization constant depending on $\nu
\mbox{ }\left( {0<\nu <1} \right)$.

Following Eq. (\ref{eq1}), Losada and Nieto suggested the new fractional
Caputo-Fabrizio derivative operator [10-12]
\begin{equation}
\label{eq2}
{ }_\ast ^{CF} D_x^{\left( \nu \right)} T\left( x \right)=\frac{1}{1-\nu
}\int\limits_0^x {\exp \left( {-\frac{\nu }{1-\nu }\left( {x-\lambda }
\right)} \right)T^{\left( 1 \right)}\left( \lambda \right)d\lambda } ,
\end{equation}
where $\nu \mbox{ }\left( {0<\nu <1} \right)$ is a real number and $\Im
\left( \nu \right)=2/\left( {2-\nu } \right)$.

Eqs. (\ref{eq1}) and (\ref{eq2}) represent an extension of the Caputo fractional derivative
with singular kernel. However, an analog of the Riemann-Liouville fractional
derivative with singular kernel has not yet been formulated. The main aim of
the article is to propose a new fractional derivative without singular
kernel, which is an extension of the Riemann-Liouville fractional derivative
with singular kernel, and to study its application in the modeling of the
fractional-order heat flow.

In this line of thought, the structure of the article is as follows. Section
2 presents a new fractional derivative without singular kernel. Section 3
discusses its application to the steady fractional-order steady heat flow in
the heat-conduction problem. Finally, Section 4 outlines the conclusions.

\textbf{Mathematical tools }

The Riemann-Liouville fractional derivative of fractional order $\nu $ of
the function $T\left( x \right)$ is defined as [1]
\begin{equation}
\label{eq3}
^{RL}D_{a^+}^{\left( \nu \right)} {\rm T}\left( x \right)=\frac{1}{\Gamma
\left( {1-\nu } \right)}\frac{d}{dx}\int\limits_a^x {\frac{{\rm T}\left(
\lambda \right)}{\left( {x-\lambda } \right)^\nu }d\lambda } ,
\end{equation}
where $a\le x$ and $\nu \mbox{ }\left( {0<\nu <1} \right)$ is a real number.

Replacing the function $1/\left( {x-\lambda } \right)^\nu \Gamma \left(
{1+\nu } \right)$ by the function
\[\Re \left( \nu \right)\exp \left(
{-\frac{\nu }{1-\nu }\left( {x-\lambda } \right)} \right)\]

$/\left( {1-\nu } \right),$ we obtain a new fractional derivative given by:
\begin{equation}
\label{eq4}
D_{a^+}^{\left( \nu \right)} {\rm T}\left( x \right)=\frac{\Re \left( \nu
\right)}{1-\nu }\frac{d}{dx}\int\limits_a^x {\exp \left( {-\frac{\nu }{1-\nu
}\left( {x-\lambda } \right)} \right){\rm T}\left( \lambda \right)d\lambda }
,
\end{equation}
where $a\le x$, $\nu \mbox{ }\left( {0<\nu <1} \right)$ is a real number and
$\Re \left( \nu \right)$ is a normalization function depending on $\nu $
such that $\Re \left( 0 \right)=\Re \left( 1 \right)=1$.

Taking $\psi =1/\nu -1$, with $0<\psi <+\infty $, Eq. (\ref{eq4}) can be rewritten
as:
\begin{equation}
\label{eq5}
D_{a^+}^{\left( {\frac{1}{\psi +1}} \right)} {\rm T}\left( x \right)=\aleph
\left( \psi \right)\frac{d}{dx}\int\limits_a^x {\Pi \left( \lambda
\right){\rm T}\left( \lambda \right)d\lambda } ,
\end{equation}
where $\aleph \left( \psi \right)=\left( {\psi +1} \right)\Re \left(
{1/\left( {\psi +1} \right)} \right)$ and $\Pi \left( \lambda \right)=\exp
\left( {-\left( {x-\lambda } \right)/\psi } \right)/\psi $.

With the help of the~following approximation to the identity [9, 13]
\begin{equation}
\label{eq6}
\mathop {\lim }\limits_{\psi \to 0} \Pi \left( \lambda \right)=\delta \left(
{x-\lambda } \right),
\end{equation}
where $\nu \to 1$ (or $\psi \to 0)$, Eq. (\ref{eq4}) becomes
\begin{equation}
\label{eq7}
\mathop {\lim }\limits_{\nu \to 1} D_{a^+}^{\left( \nu \right)} {\rm
T}\left( x \right)=\mathop {\lim }\limits_{\psi \to 0} \aleph \left( \psi
\right)\frac{d}{dx}\int\limits_a^x {\Pi \left( \lambda \right){\rm T}\left(
\lambda \right)d\lambda } ={\rm T}^{\left( 1 \right)}\left( x \right).
\end{equation}
When $\nu \to 0$ (or $\psi \to +\infty )$, Eq.(\ref{eq4}) can be written as
\begin{equation}
\label{eq8}
\mathop {\lim }\limits_{\nu \to 0} D_{a^+}^{\left( \nu \right)} {\rm
T}\left( x \right)=\mathop {\lim }\limits_{\nu \to 0} \frac{\Re \left( \nu
\right)}{1-\nu }\frac{d}{dx}\int\limits_a^x {\exp \left( {-\frac{\nu }{1-\nu
}\left( {x-\lambda } \right)} \right){\rm T}\left( \lambda \right)d\lambda }
={\rm T}\left( x \right).
\end{equation}
Taking the Laplace transform of the new fractional derivative without
singular kernel for the parameter $a=0$, we have
\begin{equation}
\label{eq9}
L\left( {D_0^{\left( \nu \right)} {\rm T}\left( x \right)} \right)=\frac{\Re
\left( \nu \right)s}{\nu \left( {1-s} \right)+s}{\rm T}\left( s \right),
\end{equation}
where $L\left( {\xi \left( x \right)} \right):=\int\limits_0^x {\exp \left(
{-sx} \right)\xi \left( x \right)dx} =\xi \left( s \right)$ represents the
Laplace transform of the function $\xi \left( x \right)$ (see [14]).

We now consider
\begin{equation}
\label{eq10}
{\rm T}\left( s \right)=\left( {\frac{\nu }{\Re \left( \nu
\right)s}+\frac{1-\nu }{\Re \left( \nu \right)}} \right)\Xi \left( s
\right),
\end{equation}
where $D_0^{\left( \nu \right)} {\rm T}\left( x \right)=\Xi \left( x
\right)$ and $L\left( {\Xi \left( x \right)} \right)=\Xi \left( x \right)$.

Taking the inverse Laplace transform of Eq. (\ref{eq10}) we obtain
\begin{equation}
\label{eq11}
{\rm T}\left( x \right)=\frac{1-\nu }{\Re \left( \nu \right)}\Xi \left( x
\right)+\frac{\nu }{\Re \left( \nu \right)}\int_0^x {\Xi \left( x \right)dx}
,\mbox{ }x>0,\mbox{ }0<\nu <1.
\end{equation}
If $0<\nu <1$ and $\Re \left( \nu \right)=1$, then Eq. (\ref{eq4}) and Eq. (\ref{eq11}) can
be written as
\begin{equation}
\label{eq12}
_\ast D_{a^+}^{\left( \nu \right)} {\rm T}\left( x \right)=\frac{1}{1-\nu
}\frac{d}{dx}\int\limits_a^x {\exp \left( {-\frac{\nu }{1-\nu }\left(
{x-\lambda } \right)} \right){\rm T}\left( \lambda \right)d\lambda } ,
\end{equation}
and
\begin{equation}
\label{eq13}
{\rm T}\left( x \right)=\left( {1-\nu } \right)\Xi \left( x \right)+\nu
\int_0^x {\Xi \left( x \right)dx} ,\mbox{ }x>0,\mbox{ }0<\nu <1,
\end{equation}
respectively.

\textbf{Modelling the fractional-order steady heat flow }

The fractional-order Fourier law in one-dimension case is suggested as:
\begin{equation}
\label{eq14}
KD_0^{\left( \nu \right)} {\rm T}\left( x \right)=-H\left( x \right),
\end{equation}
where $K$ is the thermal conductivity of the material and $H\left( x
\right)$ represents the heat flux density.

The heat flow of the fractional-order heat conduction is presented as
\begin{equation}
\label{eq15}
H\left( x \right)=g,
\end{equation}
where $g$ is the heat flow (a constant) of the material.

By submitting Eq. (\ref{eq13}) into Eq. (\ref{eq14}) and taking the Laplace transform, it
results:
\begin{equation}
\label{eq16}
\frac{\Re \left( \nu \right)s}{\nu +\left( {1-\nu } \right)s}{\rm T}\left( s
\right)=-\frac{g}{K},
\end{equation}
which leads to
\begin{equation}
\label{eq17}
{\rm T}\left( s \right)=\frac{-g\left( {\nu +\left( {1-\nu } \right)s}
\right)}{K\Re \left( \nu \right)s}.
\end{equation}
Taking the inverse Laplace transform of Eq. (\ref{eq17}), we obtain
\begin{equation}
\label{eq18}
{\rm T}\left( x \right)=-C\left( {\frac{g\nu x}{K\Re \left( \nu
\right)}+\frac{g\left( {1-\nu } \right)}{K\Re \left( \nu \right)}} \right),
\end{equation}
where $C$ is a constant depending on the initial value ${\rm T}\left( x
\right)$.

The corresponding graphs with different orders $\nu =\left\{ {0.3,\mbox{
}0.6,\mbox{ }1} \right\}$ are shown in Figure 1.

\begin{figure}[htbp]
\centerline{\includegraphics[width=5.83in,height=4.38in]{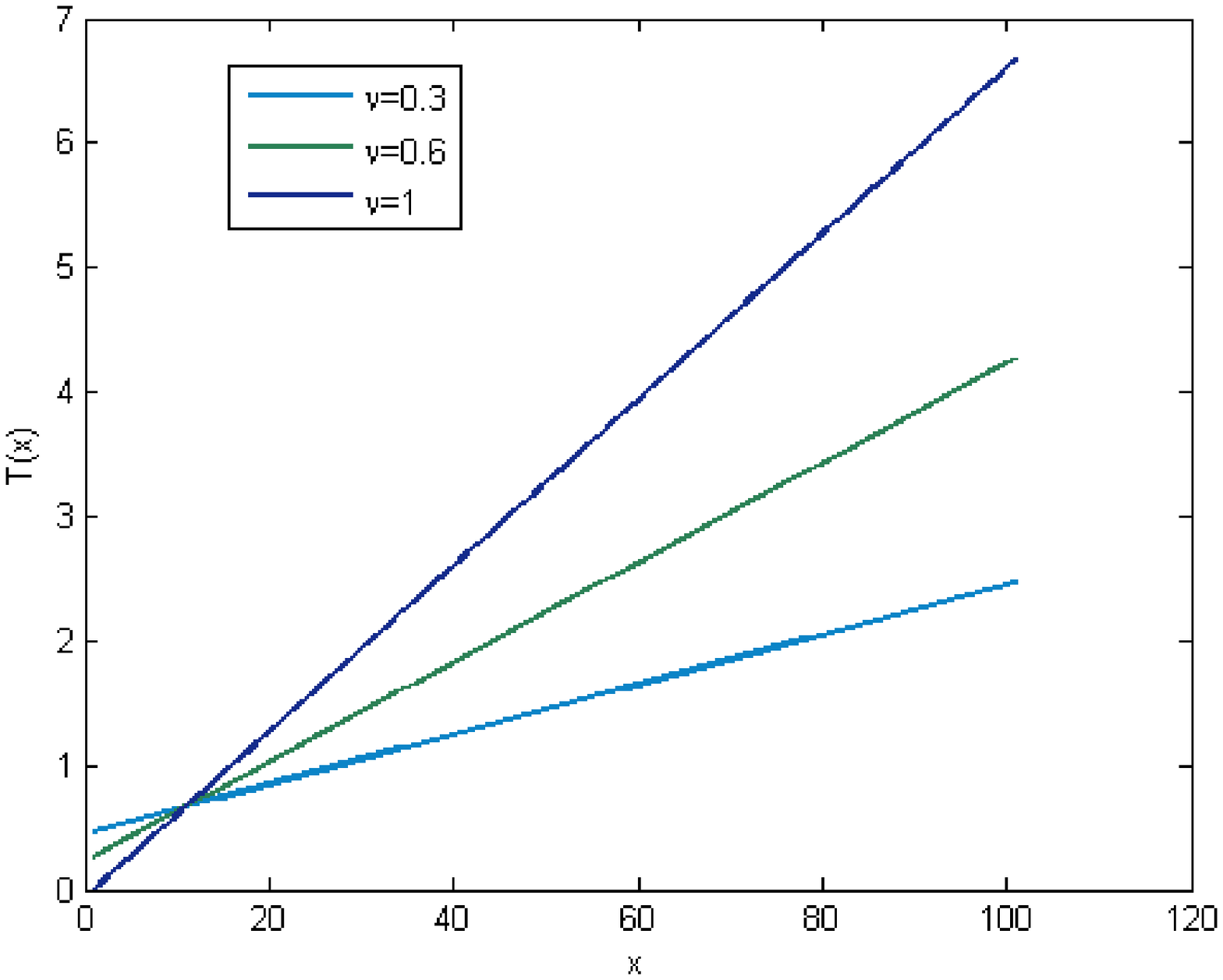}}
\label{fig1}
\end{figure}

\begin{center}
\textbf{Figure 1. The plots of }${\rm T}\left( x \right)$\textbf{ with the
parameters }$\nu =\left\{ {0.3,\mbox{ }0.6,\mbox{ }1} \right\}$\textbf{,
}$C=-1$\textbf{, }$g=2$\textbf{, }$K=3$\textbf{ and }$\Re \left( \nu
\right)=1$\textbf{. }
\end{center}

\textbf{Conclusions}

In this work a new fractional-order operator without singular kernel, which
is an analog of the Riemann-Liouville fractional derivative with singular
kernel, was proposed for the first time. An illustrative example for
modelling the fractional-order steady heat flow was given and the analytical
solution for the governing~equation~involving the fractional derivative
without singular kernel was discussed.

\textbf{References}

\begin{enumerate}
\item Yang, X. J, Baleanu, D, Srivastava, H. M., \textit{Local Fractional Integral Transforms and Their Applications}, Academic Press, New York, 2015
\item Hristov J, et al., Thermal impedance estimations by semi-derivatives and semi-integrals: 1-D semi-infinite cases, \textit{Thermal Science}, 17 (2013), 2, pp.581-589
\item Povstenko Y. Z., Fractional radial heat conduction in an infinite medium with a cylindrical cavity and associated thermal stresses, \textit{Mechanics Research Communications}, 37(2010), 4, pp.436-40
\item Hussein, E. M., Fractional order thermoelastic problem for an infinitely long solid circular cylinder, \textit{Journal of Thermal Stresses}, 38 (2015), 2, pp.133-45
\item Wei, S., et al., Implicit local radial basis function method for solving two-dimensional time fractional diffusion equations, \textit{Thermal Science}, 19 (2015), S1, pp.59-67
\item Zhao, D., et al., Some fractal heat-transfer problems with local fractional calculus, \textit{Thermal Science}, 19 (2015), 5, pp.1867-1871
\item Yang, X. J, et al., A new numerical technique for solving the local fractional diffusion equation: Two-dimensional extended differential transform approach, \textit{Applied Mathematics and Computation}, 274 (2016), pp.143-151
\item Jafari, H., et al., A decomposition method for solving diffusion equations via local fractional time derivative, \textit{Thermal Science}, 19 (2015), S1, pp.123-9
\item Caputo, M., et al., A new definition of fractional derivative without singular Kernel,~\textit{Progress in Fractional Differentiation and Applications},~1 (2015), 2, pp.73-85
\item Lozada, J. et al., Properties of a new fractional derivative without singular kernel,~\textit{Progress in Fractional Differentiation and Applications}, 2015(\ref{eq1}),~2, pp.87-92
\item Atangana, A., On the new fractional derivative and application to nonlinear Fisher's reaction--diffusion equation, \textit{Applied Mathematics and Computation}, 273(2016), pp.948-956
\item Alsaedi, A., et al., Fractional electrical circuits, \textit{Advances in Mechanical Engineering}, 7 (2015),12, pp.1-7
\item Stein, E., Weiss, G.,~\textit{Introduction to Fourier Analysis on Euclidean Spaces}, Princeton University Press, 1971
\item Debnath, L., Bhatta, D., \textit{Integral Transforms and Their Applications}, CRC Press, 2014.
\end{enumerate}

\end{document}